\documentclass{article}
\usepackage{amsfonts,amsmath}
\mathchardef\mhyphen="2D % Define a "math hyphen"
\usepackage[a4paper]{geometry}
\usepackage{microtype}

\def\zo/{$0\mkern2mu\mhyphen1$}

\def\nn/{$n \times n$}

\title{A Set of Conjectured Identities for Stirling Numbers of the First Kind}
\date{\today}
\author{Paul Federbush\\
Department of Mathematics\\
University of Michigan\\
Ann Arbor, MI, 48109-1043}
\begin{document}
\maketitle

\begin{abstract}
Given an integer $g$, $ g \geq 2$, an integer $w$, $0 \leq w \leq g - 2$, and a set of g distinct numbers, $c_1, ..., c_g$, we present a conjectured identity for Stirling numbers of the first kind. We have proven all the equalities in case $g \leq 6$; and for the case $g = 7$, provided $w \leq 3$. These expressions arise from an aspect of the study of the dimer-monomer problem on regular graphs.
\end{abstract}

We organize the paper into two parts. In the first part we present the conjectured identities . In the second part we sketch their origins in the monomer-dimer problem, arising from work of M. Pernici, \cite{per}. A brief discussion of proofs is given at the end of Part 1.\\\\
\underline{Part 1 Conjectured Identities}

The (unsigned) Stirling numbers of the first kind, $\begin{bmatrix} a \\ b \end{bmatrix}$ , are defined by

\begin{equation}
\label{1}
x(x+1)...(x+n-1) = \sum_{k=0}^{n} \begin{bmatrix} n \\ k \end{bmatrix} x^k
\end{equation}
See \cite{b}. It is easy to show $\begin{bmatrix} n \\ n-w \end{bmatrix}$ is a polynomial in n of degree 2w. So we may naturally define $\begin{bmatrix} x \\ x-w \end{bmatrix}$ for any number x by extending the domain of the polynomial. We set 
\begin{equation}
\label{2}
P_{w}(x) \equiv \begin{bmatrix} x \\ x-w \end{bmatrix}
\end{equation}

Now we give ourself an integer $g \geq 2$, an integer $w$, $0\leq w \leq g-2$, and a set of $g$ distinct numbers, 
\begin{equation}
\label{3}
S =\{ c_1, ..., c_g\}
\end{equation}
We define a \underline{configuration} as a sequence of non-empty subsets of S
\begin{equation}
\label{4}
S_1,S_2,...,S_r
\end{equation}
that are disjoint with union  $S$, i.e.
\begin{equation}
\label{5}
S_{i} \neq \emptyset, S_{i} \cap S_{j} = \emptyset \text{ if } i \neq j, \bigcup_{i=1}^{r}S_i = S
\end{equation}
For a configuration we define
\begin{equation}
\label{6}
t_i = \sum_{c_{i}\in S_{i}} c_{i}, i = 1,...,r
\end{equation}
A \underline{weighted configuration} is a configuration as above for which each $S_{i}$ is assigned a non-negative integer, $w_{i}$, its weight, with the restriction
\begin{equation}
\label{7}
\sum_{i=1}^{r} w_i = w
\end{equation}
Such a weighted configuration has an \underline{evaluation} defined as
\begin{equation}
\label{8}
(-1)^{r} \frac{1}{r} \Pi_{i} P_{w_i}(t_{i})
\end{equation}

\noindent\fbox{%
    \parbox{\textwidth}{%
        The conjectured identity is that the sum over all distinct weighted configurations of their evaluations is zero.

    }%
}

We have proved the identities in certain cases by directly computing the sum of evaluations. In this way we verified the identities for $g \leq 4$ by hand. In addition we did all cases for $g \leq 6$ by computer, and likewise all cases with $g = 7$ and $w \leq 3$. The computer computations were done on a desktop computer using Maple in a couple of days. Clearly one could go further by computer. We have no idea how to organize a general proof (if one exists), we expect any general proof to be extremely difficult. It is possible that someone may find inspiration to a line of proof from the source of the identities sketched in Part 2.\\\\
\underline{Part 2 Through the Monomer-Dimer Problem, More General Conjectured Identities}

I. Wanless developed a formalism for the monomer-dimer problem on regular graphs that is the source of all our conjectured identities, \cite{Wau}. More particularly we take our formulas from the work of M. Pernici \cite{per}, a systemization of Wanless's results. We extract from this treatment of the monomer-dimer problem only the formulas that directly lead to our identities. We feel this is the only material that might stimulate ideas towards a proof.

We use three equations, (10), (12), and (16) from \cite{per}, with slight modification:
\begin{equation}
\label{9}
M_{j}(n,r) = [x^{j}]\text{exp}(nrx- \sum_{s=2} \frac{n u_{s}(r)}{s}(-x)^{s})
\end{equation}
\begin{equation}
\label{10}
M_{j} = \frac{n^{j}r^{j}}{j!} \sum_{h=0}^{j-1} \frac{a_{h}(r,j)}{n^{h}}
\end{equation}
\begin{equation}
\label{11}
[j^{k}n^{-h}]\text{ln}(1+ \sum_{s=1}^{j-1} \frac{a_{s}(r,j)}{n^{s}}) = 0, k \geq h +2
\end{equation}
The symbol $[x^{j}]$ extracts from the expression following it the coefficient of the $x^{j}$ term in its expansion in powers of x. Likewise the symbol $[j^{k}n^{-h}]$ extracts from the expression following it the coefficient of the $j^{k}n^{-h}$ term in an expansion in powers of j and $\frac{1}{n}$.

In \cite{per} Pernici with a computation using an idealized physical argument derives equation (11) from equations (9) and (10) for particular values of the $u_{s}(r)$, $s \geq 2$. This result is checked in \cite{per} for a large number of cases. We now state a set of conjectured identities more general than those in Part 1.

\noindent\fbox{%
    \parbox{\textwidth}{%
        For any set of values $u_{s}(r), s\geq 2$ equation (11) follows from equation (9) and (10).

    }%
}

A similar set of conjectured identities is given in Section 10 of \cite{Fed},  both of these sets of identities have been checked by computer computation for many cases.

We can derive the conjectured identities of Part 1, that is the boxed statement after eq(8), from the boxed statement just above as follows. Let us treat the conjectured identity from Part 1 corresponding to parameter values $g,w$, $c_1,...,c_g$ where the $c_{i}$ are distinct integers $\geq 2$.
Upon some reflection it is easy to see that verifying the identities of this type suffices to verify all the identities of Part 1.

We set
\begin{equation}
\label{12}
k = \sum c_i - w
\end{equation}
and
\begin{equation}
\label{13}
h = \sum c_i - g
\end{equation}
so that $k - h = g - w \geq 2$. And we set all the $u_{s}$ to zero except for values of $s$ equal to one of the $c_i$. Then the left side of (11) is identically equal to zero as a function of such $u_{s}$. We look at the coefficient of the term $u_{c_{1}}u_{c_{2}}...u_{c_{g}}$ in a power series expansion of the left side of (11) in the $u_{s}$. The statement that this coefficient is zero is the desired conjectured identity of Part 1. It is perhaps not trivial to verify this. ( We think it should not be too difficult to prove the implication in the other order, that
the identities of Part 1 imply the identities of Part 2. )\\
\underline{Part 3 Postscript}

I find the three sets of conjectured identities (those of Part 1, those of Part 2, and those of \cite{Fed}) very mysterious. A colleague of mine has said of one of the sets that he thinks once a proof is found it will be easy. My feeling is the opposite, that they are candidates for being things true but not provable.

ACKNOWLEDGMENT I would like to thank John Stembridge for much useful information.

\end{document}